\documentclass[12pt]{article}

\usepackage[latin1]{inputenc}
\usepackage{amsmath,amsthm,amssymb,amscd,}
\newtheorem{thm}{Theorem}[section]
 \newtheorem{cor}[thm]{Corollary}
 \newtheorem{lem}[thm]{Lemma}
 \newtheorem{prop}[thm]{Proposition}
 \theoremstyle{definition}
 \newtheorem{df}[thm]{Definition}
 \theoremstyle{remark}
 \newtheorem{rem}[thm]{Remark}
 \newtheorem{ex}{Example}
 \numberwithin{equation}{section}

\def\be#1 {\begin{equation} \label{#1}}
\newcommand{\ee}{\end{equation}}
\def\dem {\noindent {\bf Proof : }}

\def\sqw{\hbox{\rlap{\leavevmode\raise.3ex\hbox{$\sqcap$}}$%
\sqcup$}}
\def\findem{\ifmmode\sqw\else{\ifhmode\unskip\fi\nobreak\hfil
\penalty50\hskip1em\null\nobreak\hfil\sqw
\parfillskip=0pt\finalhyphendemerits=0\endgraf}\fi}

\newcommand{\mb}{\medskip\noindent}
\newcommand{\gb}{\bigskip\noindent}
\newcommand{\R}{\mathbb R}

\newcommand{\N}{\mathbb N}
\newcommand{\Z}{\mathbb Z}

\newcommand{\s}{\mathcal S}

\newcommand{\T}{{\bf T}}
\newcommand{\ot}{{{t}}}
\newcommand{\OQ}{\mathbf Q}
\newcommand{\OP}{\mathbf P}
\newcommand{\po}{{\mathbb P}}

\title{Local estimates and global continuities in Lebesgue spaces for bilinear operators. }
\author{Fr\'ed\'eric Bernicot \\
 \small Universit\'e de Paris-Sud, Orsay et CNRS 8628, 91405 Orsay Cedex, France \\
\small {\em E-mail address:} {Frederic.Bernicot@math.u-psud.fr} }

\begin{document}

\maketitle

\begin{abstract}
 In this paper, we first prove some local estimates for bilinear operators (closely related to the bilinear Hilbert
 transform and similar singular operators) with truncated symbol. Such estimates, in accordance with the Heisenberg uncertainty principle correspond to a description of ``off-diagonal'' decay. In addition they allow us to prove
 global continuities in Lebesgue spaces for bilinear operators with spatial dependent symbol. 
\end{abstract}

\tableofcontents

\section{Introduction}

The simplest bilinear operator is the product $\Pi$ , defined by
$$ \forall\, f,g\in\s(\R), \qquad \Pi(f,g)(x):=f(x)g(x).$$
Hölder's inequalities give us the continuities on Lebesgue spaces for this operator. So for all exponents
$0<p,q,r\leq \infty$ such that

\be{homogeneite} \frac{1}{p}+\frac{1}{q}=\frac{1}{r},\ee

\mb the operator $\Pi$ is continuous from $L^p(\R) \times L^q(\R)$ into $L^r(\R)$. Also a natural question appears :
how
can we modify this bilinear operation and simultaneous keep these continuities ? \\
First let $T$ be a bilinear operator, acting from $\s(\R) \times \s(\R)$ into $\s'(\R)$. It is well-known that we have
a spatial representation of $T$ with a kernel $K\in\s'(\R^3)$ and a frequential representation with a symbol
$\sigma\in\s'(\R^3)$ such that (in distributional sense) \begin{align}
 \forall\, f,g\in\s(\R), \qquad T(f,g)(x)& =\int_{\R^2}
K(x,y,z) f(y) g(z) dydz \nonumber \\
 & = \int_{\R^2} e^{ix(\alpha+\beta)} \sigma(x,\alpha,\beta) \widehat{f}(\alpha)
\widehat{g}(\beta) d\alpha d\beta. \label{operateur} \end{align} In the rest of this paper, we denote by $T_\sigma$ the
operator associated to the symbol $\sigma$. The kernel and the symbol are related by the relation
\be{noyau-symbole} K(x,y,z) = \int_{\R^2} e^{i[\alpha(x-y)+\beta(x-z)]}\sigma(x,\alpha,\beta) d\alpha d\beta.\ee
For example, the product operator $\Pi$ is given by the symbol $\sigma(x,\alpha,\beta)=1$. \\

\mb A first class of bilinear symbols is the one satisfying the ``bilinear Hörmander's condition''~: \be{condi1} \forall\, a,b,c\geq 0, \qquad \left| \partial_x^a \partial_\alpha^b \partial_\beta^c
\sigma(x,\alpha,\beta) \right| \lesssim \left(1+|\alpha|+|\beta| \right)^{-b-c}. \ee The operators $T_\sigma$ were already studied par R. Coifman and Y. Meyer in \cite{cm,cm2}, C. Kenig and E.M. Stein in \cite{KS} and recently by L. Grafakos and R. Torres in \cite{GT}. We know that under (\ref{condi1}), the operator $T_\sigma$ is bounded from $L^{p}(\R) \times L^q(\R)$ in
$L^{r}(\R)$ for all exponents $p,q,r$ verifying (\ref{homogeneite}) and $1<p,q<\infty$. In fact if the symbol is
$x$-independent, one can just assume an homogeneous decay in (\ref{condi1}) (i.e. with $(|\alpha|+|\beta|)^{-b-c}$) and then these operators can be decomposed
with paraproducts, which were first exploited by J.M. Bony in \cite{bony} and R. Coifman and Y. Meyer in \cite{cm}. The
paraproducts are studied with the linear tools (the Calder\'on-Zygmund decomposition, the Littlewood-Paley theory and
the concept of Carleson measure). In order to get the continuities for $x$-dependent symbols, pointwise estimates of
the bilinear kernel are used. Mainly for a symbol $\sigma$ satisfying (\ref{condi1}), integrations by parts allow us to
obtain \be{pointwisekernel} \left|K(x,y,z)\right| \lesssim (1+|x-y|+|x-z|)^{-M} \ee for any large enough integer $M$.
This estimate is very useful and precisely describes the ``off-diagonal'' decrease of the operator. Such an information
helps us to reduce the study of $x$-dependent symbols to the study of $x$-independent symbols (and so to the study of
paraproducts). Through these ideas, this first class of symbols are well understood nowadays. We note that this
reduction (using pointwise estimates on the kernel) has already been used in the linear case
to study the pseudo-differential operators of the well-known class $op(S^{0}_{1,0})$. Thus ``off-diagonal'' estimates play an important role.

\mb Since the work of A.Calder\'on about the $L^2$ boundedness of commutators and Cauchy integrals
(\cite{calderon1,calderon2} in the 70's), more singular bilinear operators have appeared. Mainly, he shew that the
commutators and Cauchy integrals can been decomposed by using the bilinear Hilbert transforms. The bilinear Hilbert
transform $H_{\lambda_1,\lambda_2}$ is defined by
$$ \forall\, f,g\in\s(\R), \qquad  H_{\lambda_1,\lambda_2}(f,g)(x):= p.v. \int_\R f(x-\lambda_1y)g(x-\lambda_2y)\frac{dy}{y}.$$
The $x$-independent symbol is 
$$ \sigma(\alpha,\beta)=i\pi sign(\lambda_1\alpha+\lambda_2\beta) $$ 
and so is singular on a whole line in the frequency plane. 
A.Calder\'on  conjectured that these operators are continuous on Lebesgue spaces. This famous conjecture was first
partially solved by M. Lacey and C. Thiele in \cite{LT2,LT1,LT3,LT4}. Then L. Grafakos and X. Li show in\cite{GL,Li}
some uniform (with respect to the parameters $\lambda_1$ and $\lambda_2$) continuities. These proofs use a technical time
frequency analysis, which was proven by C. Muscalu, T. Tao and C. Thiele in \cite{MTT2,mtt,MTT3} and independently by
J. Gilbert and A. Nahmod in \cite{GN,GN2}. They also get a very important result in the study of bilinear operators :
continuities in Lebesgue spaces for more singular operators than those of the first class. We are interested by these
bilinear operators and we will deal with them and some ``smooth spatial perturbations''. So we replace in
(\ref{condi1}) the quantity $|\alpha|+|\beta|=d((\alpha,\beta),0)$ by the lower quantity $d((\alpha,\beta),\Delta)$,
where $\Delta$ is a line in the frequency plane~:
$$ \Delta:= \left\{ (\alpha,\beta)\in\R^2, \lambda_1\alpha + \lambda_2\beta=0 \right\}.$$
We assume that $\Delta$ is nondegenerate, i.e. $\lambda_1$ and $\lambda_2$ are non vanishing reals and not equal, in order that $\Delta$ be a graph over the
three variables $\alpha$, $\beta$ and $\alpha+\beta$. We assume that the symbol $\sigma$ satisfies \be{condi2}
\forall\, a,b,c\geq 0, \qquad \left|
\partial_x^a
\partial_\alpha^b
\partial_\beta^c \sigma(x,\alpha,\beta) \right| \lesssim \left(1+|\lambda_1\alpha+\lambda_2\beta| \right)^{-b-c}. \ee
In the previous mentioned papers, the main result is the following one~: if $\sigma$ is $x$-independent and satisfies
(\ref{condi2}) (or the homogeneous version) then $T_\sigma$ is continuous from $L^p(\R) \times L^q(\R)$ in $L^r(\R)$
for every exponents $0<p,q,r\leq\infty$ verifying
$$ 0<\frac{1}{r}=\frac{1}{p}+\frac{1}{q}<\frac{3}{2}\quad \textrm{and} \quad 1<p,q\leq \infty.$$
So there is a natural question (asked in \cite{bipseudo}) : if an $x$-dependent symbol satisfies (\ref{condi2}), is the operator $T_\sigma$ continuous from $L^p(\R) \times L^q(\R)$ into $L^r(\R)$ with the same exponents $p,q,r$ ? In \cite{t1bilineaire}, A. Benyi, C. Demeter, A. Nahmod, R. Torres, C. Thiele and P. Villarroya proved a general result for singular integrals kernels. As an example, they can apply their result to pseudo-differential
operators associated to symbols $\sigma(x,\alpha,\beta)=\tau(x, \lambda_1 \alpha + \lambda_2 \beta)$ with $\tau$ in the class $S^0_{1,0}$ because of a modulation invariant condition imposed. Here we are able to treat general symbols satisfying (\ref{condi2}) and complete the answer to the question in \cite{bipseudo}. These operators do not fall under the scope of \cite{t1bilineaire} because they do not have modulation invariance. On the other hand, the general operators in \cite{t1bilineaire} cannot be realized as pseudo-differential bilinear operators with symbols
verifying (\ref{condi2}) because of the minimal regularity assumptions required in the kernels. \\

\mb With this aim, we would like to use the same arguments as for the symbols satisfying (\ref{condi1}), where we have
seen the important role of the ``off-diagonal'' decay of the bilinear kernel, obtained with integrations by
parts. For our more singular operators, integration by parts does not work : to obtain a description of
``off-diagonal'' estimates  is the most important difficulty.

\mb We now come to our main result. For notation, we denote the norm in $L^p(E)$ for any measurable set $E\subset \R$ by $\|\ \|_{p,E,dx}$ (or $\|\ \|_{p,E}$
 if there is no confusion for the variable). For an interval $I$, we set
\be{couronne} C_k(I):= \left\{ x\in \R,\ 2^{k} \leq 1+\frac{d(x,I)}{|I|} <2^{k+1} \right\} \ee the scaled corona around
$I$. So we have $C_0(I)=2I$ and $C_{k}(I) \subset 2^{k+1}I$. We will first prove the following result~:

\begin{thm}  \label{thm:central} Let $\Delta$ be a nondegenerate line of the frequency plane. Let $p,q$ be exponents such that
$$1< p,q \leq \infty \quad \textrm{and} \quad  0<\frac{1}{r}=\frac{1}{q}+\frac{1}{p} <\frac{3}{2}.$$
Then for all $\delta\geq 1$, there is a constant $C=C(p,q,r,\Delta,\delta)$  such that for all interval $I\subset \R$,
for all symbol $\sigma\in C^\infty(\R^3)$ satisfying \be{heisenberg} \forall\, a,b,c\geq 0 \qquad
\left|\partial_x^a\partial^{b}_\alpha \partial^c_\beta \sigma(x,\alpha,\beta)\right| \lesssim
\left(|I|^{-1}+d((\alpha,\beta),\Delta)\right)^{-b-c}, \ee we have the following local estimate~: for all functions
$f,g\in\s(\R)$
\begin{align*}
\left( \frac{1}{|I|} \int_{I} \left|T_\sigma(f,g)(x)\right|^r dx \right)^{1/r}  \leq C & \left[\sum_{k\geq 0} 2^{-k\delta} \left( \frac{1}{|2^{k+1}I|} \int_{C_k(I)} |f(x)|^p dx \right)^{1/p} \right] \nonumber \\
 &  \left[\sum_{k\geq 0} 2^{-k\delta} \left( \frac{1}{|2^{k+1}I|} \int_{C_k(I)} |g(x)|^{q} dx \right)^{1/q} \right].
 \end{align*}
In particular, with the Hardy-Littlewood's operator $M_{HL}$, we have
\begin{align*}
\left(|I|^{-1} \int_I |T_{\sigma} (f,g)|^r\right)^{1/r} \lesssim \inf_{I}M_{HL}\left(|f|^p\right)^{1/p}
\inf_{I}M_{HL}\left(|g|^q\right)^{1/q} \lesssim \|f\|_\infty \|g\|_\infty.
\end{align*}
\end{thm}

\mb The weight $(|I|^{-1}+d((\alpha,\beta),\Delta))^{-N}$ is not integrable over the whole frequency plane (even if
$N$ is large enough due to the modulation invariance) and therefore we cannot have a pointwise estimate of the bilinear
kernel (such (\ref{pointwisekernel}) when we assume (\ref{condi1})). So such a result is interesting because it
precisely describes ``off-diagonal'' estimates for the bilinear operator~:

\begin{cor} \label{corrolaire} With the same notations as Theorem \ref{thm:central}, for all large enough $\delta$, there exists a
constant $C=C(p,q,r,\Delta,\delta)$ such that for any measurable sets $E,F \subset \R$ we have for all functions $f\in
L^p(E)$ and $g\in L^q(F)$~: \be{decrease} \|T_\sigma(f,g)\|_{r,I} \leq C\left(1+ \frac{d(I,E)}{|I|}\right)^{-\delta}
\left(1+ \frac{d(I,F)}{|I|}\right)^{-\delta} \|f\|_{p,E}\|g\|_{q,F}. \ee
\end{cor}

\mb This corollary is a direct application of the previous Theorem. So in spite of the fact that the symbol could be
much more singular than those satisfying only (\ref{condi1}), we almost obtain the pointwise estimate
(\ref{pointwisekernel}). Here we have a description of the same fast decrease for the bilinear kernel, not with a
pointwise estimate, but with local estimates at the scale $|I|$. These local estimates are less precise than the
pointwise estimate but we will see that they are sufficient and they can have the same role.

\mb We note that Theorem \ref{thm:central} is in accordance with the Heisenberg's uncertainty principle, which tell us that if we want to localize at
the scale $|I|$ in the spatial domain, we cannot localize in the frequency domain at a lower scale than $|I|^{-1}$.
For example, ou Theorem \ref{thm:central} applies if the symbol is supported in the domain $\{(\alpha,\beta), d((\alpha,\beta),\Delta)\geq |I|^{-1} \}$ and it is this case that we consider first in the proof. In fact in (\ref{heisenberg}), we allow instead a nice behavior around the line $\Delta$. With this point of view, we could call Theorem \ref{thm:central}  : an ``high frequency estimate''. In this expression,
the term "frequency" corresponds to the distance between the point $(\alpha,\beta)$ to the line of singularity
$\Delta$. We prefer the expression ``local estimates'', because we will use the spatial fast decay in
order to get the following result.

\begin{thm} \label{thmcentral} Let $\Delta$ be an nondegenerate line of the frequency plane.
Let $p,q$ be exponents such that
$$1< p,q \leq \infty \quad \textrm{and} \quad  0<\frac{1}{r}=\frac{1}{q}+\frac{1}{p} <\frac{3}{2}.$$
For all symbol $\sigma\in C^\infty(\R^3)$ satisfying
$$ \forall\, a,b,c\geq 0 \qquad \left|\partial_x^a\partial^{b}_\alpha \partial^c_\beta \sigma(x,\alpha,\beta)\right| \lesssim \left(1+d((\alpha,\beta),\Delta)\right)^{-b-c}, $$
the associated operator $T_\sigma$ is bounded from $L^p(\R) \times L^q(\R)$ into $L^r(\R)$.
\end{thm}

\mb This result answers a question of \cite{bipseudo}. In addition it
will allow us to define a bilinear pseudodifferential calculus, based on these operators : in a next paper
\cite{pseudo}, we will define classes for bilinear pseudo-differential operators of order $(m_1,m_2)$ and study their
action on Sobolev spaces. In order to carry on the work of \cite{bipseudo}, we will give rules of symbolic calculus for
the duality and the composition and also complete the construction of a bilinear pseudo-differential calculus.

\begin{rem}
The proof of Theorem \ref{thm:central} is a shake between a localization argument and the ``classical'' time-frequency
analysis used for these bilinear operators. So it is quite easy to obtain a version of our Theorems \ref{thm:central} and \ref{thmcentral} for $(n-1)$-linear operators $T_\sigma$ with a
nondegenerate space $\Delta$ of dimension $k<n/2$, by following the ideas of \cite{MTT2}. By using the results of
\cite{terwilleger}, we are able to obtain the same results for a  multidimensional problem ; and by using the
uniform estimates of \cite{mtt}, it seems possible to obtain uniform (with respect to the nondegenerate line $\Delta$)
local estimates.
\end{rem}

\mb The plan of this article is as follows. We first prove Theorem \ref{thm:central} in Section \ref{section1} for $x$-independent symbols.
 Then in Section \ref{sectiont} we get the same result for maximal bilinear operators and we 
conclude the proof of Theorem \ref{thm:central} in the general case. Then in Section \ref{weight}, we use these
local estimates to obtain global continuities for bilinear operators in weighted Lebesgue and Sobolev spaces and in
particular we prove Theorem \ref{thmcentral}.

\section{Proof of Theorem \ref{thm:central} for $x$-independent symbol.}
\label{section1}

\gb In this section, we assume that the symbol $\sigma$ is $x$-independent and is supported on the domain~:
$$ \left\{(\alpha,\beta),\ d((\alpha,\beta),\Delta)\geq |I|^{-1} \right\}.$$
We divide the proof in two subsections. First, we will quickly recall the decomposition of our bilinear operator
$T_\sigma$ by combinatorial model sums. So we will reduce the problem to a study of the ``restricted weak type'' for
some localized trilinear forms. We will study them in the second subsection.

\subsection{Reduction to the study of discrete models.} \label{discretis}

First of all, we define and recall the time-frequency tools (see for example \cite{MTT3})~:

\begin{df} A {\it tile} is a rectangle (i.e. a product of two intervals) $I\times \omega$ of area one. A {\it tri-tile} $s$ is a
rectangle $s= I_s \times \omega_s$ of bounded area, which contains three tiles $s_{i}=I_{s_i} \times \omega_{s_i}$ for
$i=1,2,3$ such that
$$ \forall i,j\in\{1,2,3\}, \qquad I_{s_i}=I_s \quad \textrm{and} \quad i\neq j \Longrightarrow \omega_{s_i} \cap \omega_{s_j} = \emptyset.$$
A set $\{I\}_{I\in {\mathcal I}}$ of real intervals is called a {\it grid} if for all $k\in \Z$
\be{grid} \sum_{\genfrac{}{}{0pt}{}{I\in {\mathcal I}}{2^{k-1}\leq |I|\leq 2^{k+1}}} {\bf 1}_{I} \lesssim {\bf 1}_{\R}, \ee
where the implicit constant is independent of $k$ and of the grid. So a grid has the same structure than the dyadic grid. \\
Let $\OQ$ be a set of tri-tiles. It is called a {\it collection} if
\begin{itemize}
 \item $\left\{ I_s,\ s\in \OQ\right\} \textrm{  is a grid,}$
 \item ${\mathcal J}:=\left\{ \omega_s,\ s\in \OQ \right\} \bigcup_{i=1}^{3} \left\{ \omega_{s_i},\ s\in \OQ \right\} \textrm{  is a grid,} $
 \item $\omega_{s_i} \subsetneq \varpi \in {\mathcal J} \Longrightarrow \forall\, j\in\{1,2,3\},\ \omega_{s_j} \subset \varpi. $
\end{itemize}
\end{df}

\mb Now we can define the wave packet for a tile.

\begin{df} Let $\Phi$ be a smooth function such that
$$\|\Phi\|_{2}=1 \quad {\textrm and} \quad \textrm{supp}(\widehat{\Phi}) \subset [-1/2,1/2].$$
For $P=I\times \omega$ a tile, we set
$$ \Phi_{P}(x):= |I|^{-1/2} \Phi\left(\frac{x-c(I)}{|I|} \right) e^{2i\pi x c(\omega)},$$
where for $U$ an interval we denote $c(U)$ its center. So $\Phi_{P}$ is normalized in the $L^2(\R)$ space, concentrated
in space around $I$ and its spectrum is exactly contained in $\omega$.
\end{df}

\mb Nowadays it is well known (see for example \cite{BG2} and \cite{BG1}) that the operator $T_\sigma$ of Theorem
\ref{thm:central} admits a decomposition~:
\begin{align}
\lefteqn{T_\sigma(f,g)(x) :=\sum_{u=(u_1,u_2,u_3) \in \Z^{3}}  } & & \nonumber \\
 & &  (1+|u|^2)^{-N}\sum_{s\in \mathbf{S}_{u}} |I_s|^{-1/2}\epsilon_{s}(u) \langle
(\tau_{u_1}\phi)_{s_1},f \rangle \langle (\tau_{u_2}\phi)_{s_2},g\rangle (\tau_{u_3}\phi)_{s_3}(x), \label{decompositionop}
\end{align}
where $\mathbf{S}_{u}$ is a collection of tri-tiles depending on $u$, $(\epsilon_{s}(u))_{s\in {\mathbf{S}_u}}$ are
bounded reals and $N$ is an integer as large as we want. We write $\tau_{v}$ for the translation operator :
$\tau_v(f)(x)= f(x-v)$. The coefficients $\epsilon_s(u)$ are uniformly bounded
with respect to the parameter $u$ and the implicit constant in (\ref{grid}) (for the definition of a grid) is bounded by the estimates of the symbol $\sigma$. \\
By using the assumption that $\sigma$ is supported in
 $\left\{(\alpha,\beta),\
|\alpha-\beta|\geq |I|^{-1}\right\}$, we have the very important property \be{remarque} \forall u\in \Z^3,\ \forall
s\in {\mathbf S}_{u}, \quad |\omega_s|\gtrsim |I|^{-1} \ (\textrm{which is equivalent to  } |I_s|\lesssim |I|). \ee

\mb So Theorem \ref{thm:central} is a consequence of the following one~:

\begin{thm} \label{thm:centralreduit}
Let $\mathbf{S}$ be a collection of tri-tiles satisfying the property (\ref{remarque}), $(\epsilon_s)_{s\in\mathbf{S}}$
bounded reals and $(\phi^i)_{i=1,2,3}$ smooth functions whose the spectrum is contained in $[-1/2,1/2]$. We denote $T_{\mathbf{S}}$ the
bilinear operator defined by
$$ T_{\mathbf{S}}(f,g)(x):= \sum_{s\in \mathbf{S}} |I_s|^{-1/2}\epsilon_{s} \langle \phi^1_{s_1},f \rangle \langle \phi^2_{s_2},g\rangle \phi^3_{s_3}(x).$$
For the exponents $(p,q,r)$ of Theorem \ref{thm:central} and for all $\delta\geq 1$, we have the local estimate
\begin{align*}
 \left(\int_I \left| T_{\mathcal{S}}(f,g) \right|^r\right)^{1/r} \lesssim  \left(\sum_{k\geq 0} 2^{-k(1/p+\delta)} \|f\|_{p,C_k(I)}\right) \left(\sum_{k\geq 0} 2^{-k(1/q+\delta)} \|g\|_{q,C_k(I)}\right).
 \end{align*}
In addition the implicit constant depends on the functions $\phi^i$ by the parameters
$$ c_M(\phi^i):= \sup_{x\in\R} \sum_{0\leq k\leq M} (1+|x|)^{M} \left| (\phi^i)^{(k)}(x) \right|$$
for $M=M(p,q,r,\delta)$ a large enough integer.
\end{thm}

\mb We have also to prove this result. In order to show it, we need to decompose the functions $f$ and $g$ around the
interval $I$. The interval $I$ being fixed, we omit it in the notation for convenience and for $i\in \N$, we set the
corona $C_i:=C_i(I)$. With the property (\ref{remarque}), we have the following decomposition
\begin{align}
T_{\mathbf{S}}(f,g) = & \sum_{k_1,k_2\geq 0} T_{\mathbf{S},0}^{k_1,k_2}(f,g) + \sum_{\genfrac{}{}{0pt}{}{k_1,k_2\geq 0}{l\leq 0}} T_{\mathbf{S},1}
^{k_1,k_2,l}(f,g) \label{decompop}
\end{align}
with
\begin{align*}
T_{\mathbf{S},0}^{k_1,k_2}(f,g)(x) & :=  \sum_{\genfrac{}{}{0pt}{}{s\in\mathbf{S}}{I_s\subset 2I} } |I_s|^{-1/2}\epsilon_{s} \langle \phi^1_{s_1},f{\bf 1}_{C_{k_1}} \rangle \langle \phi^2_{s_2},g{\bf 1}_{C_{k_2}}\rangle \phi^3_{s_3}(x)
\end{align*}
and
\begin{align*}
T_{\mathbf{S},1}^{k_1,k_2,l}(f,g)(x) & := \sum_{\genfrac{}{}{0pt}{}{s\in\mathbf{S}}{\genfrac{}{}{0pt}{}{I_s \nsubseteq 2I}{2^{l}|I|\leq |I_s| < 2^{l+1} |I|}}} |I_s|^{-1/2}\epsilon_{s} \langle \phi^1_{s_1},f{\bf 1}_{C_{k_1}} \rangle \langle \phi^2_{s_2},g{\bf 1}_{C_{k_2}}\rangle \phi^3_{s_3}(x).
\end{align*}

\mb Due to the important property (\ref{remarque}), we only have to consider tiles $s$ with $|I_s|\leq |I|$. The other
corresponding terms ($l > 0$) cannot be studied as we are going to do, according to the Heisenberg uncertainty
principle.

\mb In the next subsection, we shall prove the following theorem~:
\begin{thm} \label{thm:res} Let $(p,q,r)$ be exponents of Theorem \ref{thm:central}. The operators $T_{\mathbf{S},i}^j$ are continuous from $L^{p}(\R)\times L^q(\R)$ into $L^{r}(I)$. For convenience, we denote  $C(T_{\mathbf{S},i}^{j}):= \|T_{\mathbf{S},i}^j\|_{L^p \times L^q \to L^r}$ and we omit the exponents.
Then these continuity bounds satisfy
\begin{align*}
C(T_{\mathbf{S},0}^{k_1,k_2}) & \lesssim c_M(\phi^1)c_M(\phi^2)c_M(\phi^3) 2^{-\delta'(k_1+k_2)} \\
C(T_{\mathbf{S},1}^{k_1,k_2,l}) & \lesssim c_M(\phi^1)c_M(\phi^2)c_M(\phi^3)
2^{-\delta'(|l|+k_1+k_2)}
\end{align*}
for any large enough real $\delta'$, with an integer $M=M(p,q,r,\delta')$.
\end{thm}

\mb We claim that Theorem \ref{thm:centralreduit} is a consequence of Theorem \ref{thm:res}. \\
{\bf Proof of Theorem \ref{thm:centralreduit} :} \\
By using Theorem \ref{thm:res} and the decomposition (\ref{decompop}), we have that for all functions $f,g\in \s(\R)$ if $r\geq 1$ then
\begin{align*}
 \lefteqn{\left\| T_{\mathbf{S}}(f,g)\right\|_{r,I} \lesssim} & & \\
 & &  \sum_{k_1,k_2 \geq 0} C(T_{\mathbf{S},0}^{k_1,k_2}) \|f{\bf 1}_{C_{k_1}}\|_p
\|g{\bf 1}_{C_{k_2}}\|_q + \sum_{\genfrac{}{}{0pt}{}{k_1,k_2\geq 0}{l\leq 0}} C(T_{\mathbf{S},1}^{k_1,k_2,l}) \|f{\bf 1}_{C_{k_1}}\|_q \|g{\bf 1}_{C_{k_2}}\|_r,
\end{align*}
and if $r<1$ then
\begin{align*}
 \lefteqn{\left\| T_{\mathbf{S}}(f,g)\right\|_{r,I}^{r}  \lesssim} & & \\
 & &  \sum_{k_1,k_2 \geq 0} C(T_{\mathbf{S},0}^{k_1,k_2})^{r} \|f{\bf 1}_{C_{k_1}}\|_p^{r} \|g{\bf 1}_{C_{k_2}}\|_q^{r} + \sum_{\genfrac{}{}{0pt}{}{k_1,k_2\geq 0}{ l\leq 0} } C(T_{\mathbf{S},1}^{k_1,k_2,l})^{r} \|f{\bf 1}_{C_{k_1}}\|_p^{r} \|g{\bf 1}_{C_{k_2}}\|_q^{r}.
\end{align*}

\gb
{\bf The case $r\geq 1$.} \\
With the estimate of $C(T_{\mathbf{S},0}^{k_1,k_2})$ and $C(T_{\mathbf{S},1}^{k_1,k_2,l})$ given by Theorem \ref{thm:res}, we obtain
\begin{align*}
\lefteqn{ \left\| T_{\mathbf{S}}(f,g)\right\|_{r,I}} & & \\
 & \hspace{0.7cm} \lesssim \sum_{k_1,k_2 \geq 0} 2^{-\delta'(k_1+k_2)} \|f{\bf 1}_{C_{k_1}}\|_p
\|g{\bf 1}_{C_{k_2}}\|_q + & \\
& \hspace{5cm} \sum_{\genfrac{}{}{0pt}{}{k_1,k_2\geq 0}{l\leq 0}} 2^{-\delta'(k_1+k_2+|l|)} \|f{\bf 1}_{C_{k_1}}\|_p \|g{\bf 1}_{C_{k_2}}\|_q & \\
 & \hspace{0.7cm} \lesssim \sum_{k_1,k_2 \geq 0} 2^{-\delta'(k_1+k_2)} \|f{\bf 1}_{C_{k_1}}\|_p
\|g{\bf 1}_{C_{k_2}}\|_q. &
\end{align*}
Hence by using that $\delta'$ is as large as we want, we can conclude.

\gb
{\bf The case $r\leq 1$.} \\
We have
\begin{align*}
 \lefteqn{\left\| T_{\mathbf{S}}(f,g)\right\|_{r,I}^r} & & \\
 & \hspace{0.7cm} \lesssim  \sum_{k_1,k_2 \geq 0} 2^{-r\delta'(k_1+k_2)} \|f{\bf 1}_{C_{k_1}}\|_p^r
\|g{\bf 1}_{C_{k_2}}\|_q^r + & \\
 & \hspace{5cm} \sum_{\genfrac{}{}{0pt}{}{k_1,k_2\geq 0}{l\leq 0}} 2^{-r\delta'(k_1+k_2+|l|)} \|f{\bf 1}_{C_{k_1}}\|_p^r \|g{\bf 1}_{C_{k_2}}\|_q^r & \\
 & \lesssim \sum_{k_1,k_2 \geq 0} 2^{-r\delta'(k_1+k_2)} \|f{\bf 1}_{C_{k_1}}\|_p^r
\|g{\bf 1}_{C_{k_2}}\|_q^r. &
\end{align*}
By using H\"older's inequality and $\rho>0$ such that $p^{-1}+\rho,q^{-1}+\rho<1$, we obtain~:
\begin{align*}
\lefteqn{\left\| T_{\mathbf{S}}(f,g)\right\|_{r,I} \lesssim} & & \\
 & &  \left(\sum_{k_1\geq 0} 2^{-k_1p(\delta'-1)(\rho+1/p)} \|f{\bf 1}_{C_{k_1}}\|_p^p \right)^{1/p} \left(\sum_{k_2\geq 0} 2^{-k_2q(\delta'-1)(\rho+1/q)} \|g{\bf 1}_{C_{k_2}}\|_q^q \right)^{1/q} \\
 &  &  \lesssim \left(\sum_{k_1\geq 0} 2^{-k_1(\delta'-1)(\rho+1/p)} \|f{\bf 1}_{C_{k_1}}\|_p \right) \left(\sum_{k_2\geq 0} 2^{-k_2(\delta'-1)(\rho+1/q)} \|g{\bf 1}_{C_{k_2}}\|_q
 \right).
\end{align*}
This corresponds to the desired result (the real $\delta'$ being as large as we want). So Theorem
\ref{thm:centralreduit} is proved in the two cases. \findem

\gb We have also reduced the proof of Theorem \ref{thm:central} (for our particular symbol $\sigma$) to the proof of
Theorem \ref{thm:res}.We will prove it in the next subsection.

\subsection{Proof of Theorem \ref{thm:res}.}

\mb By using ``duality'', to prove Theorem \ref{thm:res}, we have to estimate the trilinear form defined on
$\s(\R)\times \s(\R) \times \s(\R)$ by
 \begin{align}
 \Lambda_i^j(f_1,f_2,f_3) & := \langle T_{\mathbf{S},i}^j(f_1,f_2), f_3{\bf 1}_{I} \label{trilineaireform} \rangle \\
  & = \sum_{s \in \OQ_i^j } |I_s|^{-1/2}\epsilon_{s} \langle \phi^1_{s_1},f_1{\bf 1}_{C_{k_1}} \rangle \langle \phi^2_{s_2},f_2{\bf 1}_{C_{k_2}}\rangle \langle \phi^3_{s_3}, f_3 {\bf 1}_{I}
  \rangle,\nonumber
 \end{align}
where $\OQ_i^j$ is a collection of tri-tiles, depending on $T_{\mathbf{S},i}^{j}$. \mb We need to define usual tools of
time-frequency analysis.

\begin{df} We have already defined the tri-tiles. For $j\in\{1,2,3\}$ an index and $\ot\in\mathbf{S}$ a tri-tile, a collection $\T$ of tri-tiles is called a {\it $j$-tree} with top $\ot$ if
$$\forall\, s\in\T, \qquad I_s\subset I_{\ot} \textrm{ and } \omega_{\ot_j} \subset \omega_{s_j}.$$
Then we set $I_\T:=I_\ot$ the time-interval of the tree $\T$. A collection $\T$ of tri-tiles is called a {\it tree} if there exists an index $j\in\{1,2,3\}$ such that $\T$ is a $j$-tree. For $\T$ a $j$-tree, we define the {\it size} of the function $f_j$ over this tree by
$$size_j(\T):=\left( \frac{1}{|I_{\T}|} \sum_{s\in\T}  \left| \langle f_j,\phi^{j}_{s_j} \rangle \right| ^2\right)^{1/2}.$$
For $\OQ$ a collection of tri-tiles, we define the {\it global size} by~:
$$size_j^*(\OQ)=\sup_{\genfrac{}{}{0pt}{}{\T\subset \OQ}{ \genfrac{}{}{0pt}{}{\T=k-tree}{ k\neq j}}} size_k(\T).$$
\end{df}

\mb The quantity $|I_\T|^{1/2}size_j(\T)$ corresponds to the norm of the function $f_j$ in the space $L^2$, after
having restricted it on the tree $\T$ in the time-frequency space. \\ We recall the (abstract) Proposition 6.5 of
\cite{MTT3} (where we use Lemma 6.7 of \cite{MTT3} to estimate the quantities $\widetilde{energy_j}$)~:

\begin{prop} \label{prop} Let $(\theta_j)_{1\leq j\leq 3}$ be three exponents of $[0,1[$ satisfying
$$ \theta_1 + \theta_2 + \theta_3=1.$$
Then there exists a constant $C=C(\theta_i)$ such that for all collection $\OQ$ of tri-tiles, we have~:
$$ \left| \sum_{s \in \OQ} \frac{1}{|I_s|^{1/2}} \prod_{i=1}^{3} \langle \phi^i_{s_i},f_i \rangle \right| \leq C \prod_{i=1}^{3} size_i^*(\OQ) ^{\theta_i}   \|f_i\|_{2}^{1-\theta_i}.$$
\end{prop}

\mb This result is the main idea of this time-frequency analysis. To prove it, we use a stopping-time argument in order
to build an "orthogonal" covering of the time-frequency space with trees of ${\mathbf Q}$.

\mb Now we recall the notion of {\it restricted weak type} for trilinear forms~:

\begin{df} \label{restricteddef} For $E$ a borelian set of $\R$, we write~:
$$F(E) := \left\{ f\in\s(\R),\ \forall x\in \R,\ |f(x)|\leq {\bf 1}_{E}(x) \right\}.$$
Let $\Lambda$ be a trilinear form, defined on $\s(\R) \times \s(\R) \times \s(\R)$. Let $p_1,p_2,p_3$ be exponents of
$\R^{*}$, possiblynegative. We say that $\Lambda$ is {\it of restricted weak type} $(p_1,p_2,p_3)$ if there exists
a constant $C$ such that for all measurable sets $E_1,E_2,E_3$ of finite measure, we can find a substantial subset
$E_\beta' \subset E_\beta$ (i.e. $|E_\beta'|\geq |E_\beta|/2$) for $\beta\in\{1,2,3\}$ with \be{typerestreint}
\forall\, f_\beta \in F(E_\beta'), \qquad \left| \Lambda(f_1,f_2,f_3) \right| \leq C \prod_{\beta=1}^{3}
|E_\beta|^{1/p_\beta} \ee and $E_\beta'=E_\beta$ if $p_\beta >0$. The best constant in (\ref{typerestreint}) is called
the bound of restricted type and will be denoted by $C(\Lambda)$.
\end{df}

\mb By the real interpolation theory for trilinear forms of restricted weak type (Lemmas 3.6, 3.7, 3.8, 3.9, 3.10 and
3.11 of \cite{mtt}), Theorem \ref{thm:res} is a consequence of the following result (which is a stronger continuity
result)~:

\begin{thm} \label{thm:res2}
Let $p_1,p_2,p_3$ be non vanishing reals such that
$$ \frac{1}{p_1} + \frac{1}{p_2} + \frac{1}{p_3}=1$$
and there exists a unique index $\alpha \in \{1,2,3\}$ with $-1/2<p_\alpha^{-1}<0$ and $1/2<p_\beta^{-1}<1$ for
$\beta\neq \alpha$. Then the trilinear forms $\Lambda_i^j$ (defined by (\ref{trilineaireform})) are of restricted weak
type $(p_1,p_2,p_3)$. In addition the bounds of restricted type $C(\Lambda_i^j)$ satisfy
\begin{align*}
C(\Lambda_0^{k_1,k_2}) & \lesssim c_M(\phi^1)c_M(\phi^2)c_M(\phi^3) 2^{-\delta'(k_1+k_2)} \\
C(\Lambda_1^{k_1,k_2,l}) & \lesssim c_M(\phi^1)c_M(\phi^2)c_M(\phi^3)
2^{-\delta'(|l|+k_1+k_2)}
\end{align*}
for any real $\delta'\geq 1$ with $M=M(\delta',p_i)$ a large enough integer.
\end{thm}

\dem The exponents $(p_\beta)_\beta$ and the index $\alpha\in\{1,2,3\}$ are fixed for the proof. Let $E_1,E_2$ and $E_3$ measurable sets of finite measure. First we construct the substantial subset $E_\alpha' \subset E_\alpha$. Denote
$$ U:=\bigcup_{i=1}^3 \left\{ x\in \R,\ M_{HL}({\bf 1}_{E_i})(x)>\eta\frac{|E_i|}{|E_\alpha|}\right\}.$$
By using Hardy-Littlewood's Theorem, there exists a numerical constant $\eta$ such that
$$ |U|\leq |E_\alpha|/2.$$
We set also $E_\alpha'=E_\alpha \setminus U$. It is interesting to note that the set $E_\alpha'$ does not depend on
the form $\Lambda_i^j$. Now we fix the functions $f_\beta \in  F(E_\beta')$ for $\beta\in\{1,2,3\}$ and we shall prove
the inequality (\ref{typerestreint}). The proof is divided in three parts~: in the first step we use general estimates
for collections of tri-tiles, in the second step we will use specific estimates adapted to the above collections of
tri-tiles and then we will conclude in the third step.

\mb{\bf First step :} A general estimate. \\
Let  $\OP$ an ``abstract'' collection of tri-tiles, then for $k\geq 0$ we set $\OP_k$ the following sub-collection
$$\OP_k:=\left\{ s\in\OP, \ \ 2^{k}\leq 1+\frac{d(I_s,U^c)}{|I_s|} < 2^{k+1}\right\}. $$
These collections form a partition of $\OP$ : $\OP=\bigsqcup_{k\geq 0} \OP_k$.For each $k\geq 0$, we can apply Proposition
\ref{prop} to the collection $\OQ=\OP_k$. So for any choice of exponents $0< \theta_1,\theta_2,\theta_3<1$ with
$$ \sum_{\beta=1}^3 \theta_\beta=1, $$ we obtain
$$ \Lambda(P_k):=\left| \sum_{s\in\OP_k} |I_s|^{-1/2}\epsilon_s \prod_{\beta=1}^3 \langle f_\beta,\phi^{\beta}_{s_\beta} \rangle \right| \lesssim  \prod_{\beta=1}^3 (size_\beta^*(\OP_k))^{\theta_\beta} \|f_\beta\|_2^{1-\theta_\beta}.$$

\gb In order to estimate the quantities $size_\beta^*(\OP_k)$, we recall Lemma 7.8 of \cite{mtt}~:

\begin{lem} \label{arbre} For all integer $N$ as large as we want, there exists a constant $C=C(N)$ such that for all collection $\OQ$ of tri-tiles, for all $\beta\in\{1,2,3\}$, we have~:
$$size_\beta^*(\OQ) \leq C\sup_{s\in\OQ} \frac{1}{|I_s|} \int_\R \left(1+\frac{d(x,I_s)}{|I_s|}\right)^{-N} |f_\beta(x)| dx.$$
\end{lem}

\mb Then for $\OQ=\OP_k$, by using the definition of the sets $U$ and $E_\alpha'$, we have ~: \be{size1} \forall
\beta\neq\alpha \qquad  size_\beta^*(\OP_k) \lesssim 2^k \frac{|E_\beta|}{|E_\alpha|} \qquad \textrm{and} \qquad
size_\alpha^*(\OP_k)\lesssim 2^{-Nk}.\ee As $f_\beta$ belongs to $F(E_\beta)$, we have $\|f_\beta\|_2\leq
|E_\beta|^{1/2}$. So for $0<\epsilon<1$ and $N$ an integer as large as we want, we get~:
\begin{align*}
\Lambda(P_k) & \lesssim \prod_{\beta\neq \alpha} \left(2^k \frac{|E_\beta|}{|E_\alpha|}\right)^{\theta_\beta(1-\epsilon)}|E_\beta|^{(1-\theta_\beta)/2}  2^{Nk\theta_\alpha(1-\epsilon)} |E_\alpha|^{(1-\theta_\alpha)/2}  \prod_{\beta=1}^3 (size_\beta^*(\OP_k))^{\theta_\beta \epsilon}  \\
   & \lesssim  2^{-k} \left[ \prod_{\beta\neq\alpha} |E_\beta|^{(1+\theta_\beta)/2-\epsilon\theta_\beta}  |E_\alpha|^{(\theta_\alpha-1)/2+\epsilon(1-\theta_\alpha)}\right]  \left[\prod_{\beta=1}^3 (size_\beta^*(\OP_k))^{\theta_\beta \epsilon} \right].
 \end{align*}

\mb By definition of $size_\beta^{*}$, $\OP_k$ is a sub-collection of $\OP$ so
$$ \forall \beta\in\{1,2,3\}, \qquad size_\beta^*(\OP_k) \leq size_\beta^*(\OP).$$
We can also compute the sum over $k\geq 0$ and we obtain
\begin{align}
 \Lambda(P) & :=\left| \sum_{s\in\OP} |I_s|^{-1/2}\epsilon_s \prod_{\beta=1}^3 \langle f_\beta,\phi^{\beta}_{s_\beta} \rangle \right| \leq \sum_{k\geq 0} P_k \nonumber \\
  & \lesssim  \left[ \prod_{\beta\neq\alpha} |E_\beta|^{(1+\theta_\beta)/2-\epsilon\theta_\beta}  |E_\alpha|^{(\theta_\alpha-1)/2+\epsilon(1-\theta_\alpha)}\right] \left[\prod_{\beta=1}^3 (size_\beta^*(\OP))^{\theta_\beta \epsilon} \right]. \label{res2inter}
 \end{align}
 \mb The first term is ``good'', according to the wished global continuity. In the next step, we will use an other estimate of the quantities $size^*_\beta$,
 which will be adapted to our specific trilinear forms $\Lambda_i^j$ and which allow us to obtain the desired decays. \\

\gb {\bf Second step :} Use of the specificity of our trilinear forms $\Lambda_{i}^j$.

\mb{\underline {First case :}} the forms $\Lambda_1^{j}$. \\
In this case, we use an other decomposition~:
$$ \Lambda_1^{k_1,k_2,l}(f_1,f_2,f_3) \leq \sum_{\genfrac{}{}{0pt}{}{I_0 \nsubseteq 2I}{  2^{l-1}|I|\leq |I_0| \leq 2^{l+1}|I| }} \Lambda_1^{k_1,k_2,l}(I_0)(f_1,f_2,f_3)$$
where $I_0$ is an interval of $\R$ and
 $$ \Lambda_1^{k_1,k_2,l}(I_0)(f_1,f_2,f_3):=  \sum_{\genfrac{}{}{0pt}{}{s\in\mathbf{S}}{ I_s=I_0 } } |I_s|^{-1/2} \epsilon_s\langle f_1{\bf 1}_{C_{k_1}}, \phi^1_{s_1} \rangle \langle f_2{\bf 1}_{C_{k_2}},\phi^2_{s_2}\rangle \langle {\bf 1}_{I}f_3,\phi^3_{s_3} \rangle.$$
Let $I_0$ be fixed and denote $2^{l}=|I_0|/|I|$. The collection of tri-tiles associated to $\Lambda_1^{k_1,k_2,l}(I_0)$
is also
$$\OP:=\{ s\in \mathbf{S},\ \ I_s=I_0\}$$
For all $s\in\OP$, (by using $f_3\in F(E_3')$) we have \begin{align*} \frac{1}{|I_s|} \int_I
\left|f_3(x)\right|\left(1+\frac{d(x,I_s)}{|I_s|}\right)^{-N} dx & \leq \frac{1}{|I_s|} \int_I
\left(1+\frac{d(x,I_s)}{|I_s|}\right)^{-N} dx \\ & \leq \frac{|I|}{|I_0|} \left(1+\frac{d(I,I_0)}{|I_0|}\right)^{-N}.
\end{align*}
Then Lemma \ref{arbre} gives us
$$size_3^*(\OP) \lesssim 2^{-l} \left(1+\frac{d(I,I_0)}{|I_0|}\right)^{-N}. $$
By the same reasoning, we obtain for $f_1\in F(E_1')$ and $s\in \OP$
\begin{align*}
\frac{1}{|I_s|} \int_{C_{k_1}} \left|f_1(x)\right|\left(1+\frac{d(x,I_0)}{|I_s|}\right)^{-N} dx  & \leq \frac{1}{|I_0|} \int_{C_{k_1}} \left(1+\frac{d(x,I_0)}{|I_0|}\right)^{-N} dx \\
& \leq 2^{k_1-l} \left(1+\frac{d(C_{k_1},I_0)}{|I_0|}\right)^{-N}.
\end{align*}
And so we get
$$size_1^*(\OP) \lesssim 2^{k_1-l} \left(1+\frac{d(C_{k_1},I_0)}{|I_0|}\right)^{-N}. $$
Likely, we have
$$size_2^*(\OP) \lesssim 2^{k_2-l} \left(1+\frac{d(C_{k_2},I_0)}{|I_0|}\right)^{-N}. $$
With $\theta_1+\theta_2+\theta_3=1$ and Lemma \ref{arbre}, we can estimate~: \be{sizeA}
size_1^*(\OP)^{\theta_1}size_2^*(\OP)^{\theta_2} size_3^*(\OP)^{\theta_3} \lesssim 2^{\theta_1 k_1+ \theta_2 k_2-l}
A(I_0), \ee where $A(I_0)$ is the product of three terms
\begin{align*}
\lefteqn{A(I_0) :=} & & \\
& & \left(1+\frac{d(I,I_0)}{|I_0|}\right)^{-N\theta_3}
 \left(1+\frac{d(C_{k_1},I_0)}{|I_0|}\right)^{-N\theta_1}\left(1+\frac{d(C_{k_2},I_0)}{|I_0|}\right)^{-N\theta_2}.
 \end{align*}
We are going to get four different estimates for $A(I_0)$. \\
To keep the information about the position of $I_0$, we first have \be{A1} A(I_0) \leq
\left(1+\frac{d(I,I_0)}{|I_0|}\right)^{-N\theta_3}.\ee By using $d(I,I_0)+d(C_{k_1},I_0) \gtrsim d(I,C_{k_1}) \gtrsim
2^{k_1}|I| \simeq 2^{k_1-l}|I_0|$ and the fact that $2^{l}\leq 1$, we obtain~: \be{A2} A(I_0) \lesssim
\left(1+2^{k_1-l}\right)^{-N \min\{\theta_1,\theta_3\}} \lesssim 2^{-k_1 N \min\{\theta_1,\theta_3\}}\ee and likely
\be{A3} A(I_0) \lesssim 2^{-k_2 N \min\{\theta_2,\theta_3\}}.\ee As $I_0 \nsubseteq 2I$ and $2^l\leq 1$ then
$d(I_0,I)\geq |I|$ hence
$$\left(1+\frac{d(I,I_0)}{|I_0|}\right)^{-N} \lesssim \left(\frac{|I_0|}{|I|}\right)^N.$$
So we get \be{A4} A(I_0) \lesssim \left(\frac{|I_0|}{|I|}\right)^{N\theta_3} \lesssim 2^{lN\theta_3}.\ee Taking the
geometric mean of (\ref{A1}),(\ref{A2}), (\ref{A3}) and (\ref{A4}) (with an other exponent $N$ which is as large as we
want), we obtain~: \be{sumio} A(I_0) \lesssim 2^{-(k_1+k_2+|l|)N} \left(1+\frac{d(I,I_0)}{|I_0|}\right)^{-N}. \ee With
the help of (\ref{res2inter}) and (\ref{sizeA}), we finally estimate
\begin{align*}
\left|\Lambda_1^{k_1,k_2,l}(f_1,f_2,f_3)\right| & \lesssim \sum_{I_0} \left|\Lambda_1^{k_1,k_2,l}(I_0)(f_1,f_2,f_3)\right| \\
& \hspace{-1cm} \lesssim \sum_{I_0} \left[ \prod_{\beta\neq\alpha} |E_\beta|^{(1+\theta_\beta)/2-\epsilon\theta_\beta}  |E_\alpha|^{(\theta_\alpha-1)/2+\epsilon(1-\theta_\alpha)}\right]  2^{\epsilon(k_1+k_2+|l|)} A(I_0).
\end{align*}
The sum over the interval $I_0$ is bounded with $|I_0|=2^{l}|I|$ by (\ref{sumio}). For $N$ a large enough exponent (not
exactly the same), we have
\begin{align*}
\left|\Lambda_1^{k_1,k_2,l}(f_1,f_2,f_3) \right| \lesssim \left[ \prod_{\beta\neq\alpha}
|E_\beta|^{(1+\theta_\beta)/2-\epsilon\theta_\beta}  |E_\alpha|^{(\theta_\alpha-1)/2+\epsilon(1-\theta_\alpha)}\right]
 \tilde{C}(\Lambda_1^{k_1,k_2,l}),
\end{align*}
with
\begin{align}
 \tilde{C}(\Lambda_1^{k_1,k_2,l}) := 2^{-N\epsilon(k_1+k_2+|l|)}. \label{C}
 \end{align}

\gb {\underline{Second case :}} the forms $\Lambda_0^j$.\\
We use the same principle. We are interested by
$$ \Lambda_0^{k_1,k_2}(f_1,f_2,f_3):= \sum_{\genfrac{}{}{0pt}{}{s\in\mathbf{S}}{I_s\subset 2I }} |I_s|^{-1/2} \epsilon_{s} \langle f_1{\bf 1}_{C_{k_1}},\phi^1_{s_1} \rangle \langle f_2{\bf 1}_{C_{k_2}},\phi^2_{s_2} \rangle \langle f_3, \phi^3_{s_3} \rangle. $$
So now we choose the collection
$$\OP:=\{ s\in\mathbf{S},\ \ I_s\subset 2I\}. $$
For all $s\in\OP$,
$$ \frac{1}{|I_s|} \int_I \left(1+\frac{d(x,I_s)}{|I_s|}\right)^{-N} dx \leq 1$$ and so with Lemma \ref{arbre} we have
$$size_3^*(\OP) \lesssim 1. $$
For $f_1$, we use that
$$\frac{1}{|I_s|} \int_{C_{k_1}} \left(1+\frac{d(x,I_s)}{|I_s|}\right)^{-N} dx \lesssim
 \left(1+\frac{d(C_{k_1},I)}{|I|}\right)^{-(N-2)}$$ to conclude that
$$size_1^*(\OP) \lesssim 2^{-k_1 (N-2)}. $$
By the same argument for $f_2$, we have
$$size_2^*(\OP) \lesssim 2^{-k_2(N-2)} . $$
In this case, we can also estimate (with $N$ an other large enough integer)
\begin{align*}
size_1^*(\OP)^{\theta_1}size_2^*(\OP)^{\theta_2} size_3^*(\OP)^{\theta_3} \leq 2^{-(k_1+k_2)N\epsilon}.
\end{align*}
\mb With (\ref{res2inter}), we finally obtain
$$\Lambda_0^{k_1,k_2}(f_1,f_2,f_3)  \lesssim  \left[ \prod_{\beta\neq\alpha} |E_\beta|^{(1+\theta_\beta)/2-\epsilon\theta_\beta}  |E_\alpha|^{(\theta_\alpha-1)/2+\epsilon(1-\theta_\alpha)}\right] \tilde{C}(\Lambda_0^{k_1,k_2}),$$
with 
\be{C2} \tilde{C}(\Lambda_0^{k_1,k_2}) := 2^{-N(k_1+k_2)\epsilon}.\ee

\gb
{\bf Third step :} End of the proof. \\
For the trilinear form $\Lambda_i^j$, we have obtain a bound $C=\tilde{C}(\Lambda_i^j)$ such that for all functions $f_\beta \in F(E_\beta')$ we have
$$ \left|\Lambda_i^j(f_1,f_2,f_3) \right| \lesssim \tilde{C}(\Lambda_i^j) \left[ \prod_{\beta\neq\alpha} |E_\beta|^{(1+\theta_\beta)/2-\epsilon\theta_\beta}  |E_\alpha|^{(\theta_\alpha-1)/2+\epsilon(1-\theta_\alpha)}\right].$$
Let $(p_\beta)_\beta$ be the exponents of Theorem \ref{thm:res2}. Then we shall show that we can find
$(\theta_1,\theta_2,\theta_3)$ and $\epsilon>0$ such that
$$\forall \beta\neq\alpha \qquad (1+\theta_\beta)/2-\epsilon\theta_\beta=\frac{1}{p_\beta}  \qquad \textrm{ and } \qquad  (\theta_\alpha-1)/2+\epsilon(1-\theta_\alpha)=\frac{1}{p_\alpha}.$$
Let $\gamma>0$ be a real satisfying
$$\forall \beta\neq\alpha \qquad \left| \frac{1}{2} -\frac{1}{p_\beta}\right| < \frac{1}{2+\gamma}.$$
This is possible because $1<p_\beta<2$ for $\beta\neq \alpha$.
We begin to choose $\theta_\alpha\in]0,1[$ such that
$$ 1>\theta_\alpha >\max\left\{\theta_\alpha^0:=\frac{p_\alpha+(2+\gamma)}{p_\alpha},0\right\}$$
  and
$$ \min\left\{\frac{-1}{2+\gamma}=\frac{1}{p_\alpha(1-\theta_\alpha^0)},\frac{1}{p_\alpha} \right\} > \frac{1}{p_\alpha(1-\theta_\alpha)}> \frac{-1}{2}.$$
This is possible because $p_\alpha$ is negative and satisfies
$$\frac{1}{p_\alpha}>-\frac{1}{2}.$$
Then we get $\epsilon$ by
$$\epsilon:=\frac{1}{2}+\frac{1}{p_\alpha(1-\theta_\alpha)}\in ]0,\frac{1}{2}[ \subset ]0,1[.$$
We now define $\theta_\beta$ for $\beta\neq\alpha$ by
$$\theta_\beta:=\frac{\frac{1}{p_\beta}-\frac{1}{2}}{\frac{1}{2}-\epsilon}.$$
We have $1<p_\beta<2$ and $0<\epsilon<1/2$, so $0<\theta_\beta$ and
$$ 0<\theta_\beta = \frac{\frac{1}{p_\beta}-\frac{1}{2}}{\frac{-1}{p_\alpha(1-\theta_\alpha)}} <
\frac{\frac{1}{2+\gamma}}{\frac{1}{2+\gamma}}=1.$$ Consequently, we have solved the system of equations for the
exponents. With this choice, we obtain
$$ \forall f_1\in F(E_1'),\ f_2\in F(E_2'),\ f_3 \in F(E_3'), \qquad \Lambda_i^j(f_1,f_2,f_3) \lesssim \tilde{C}(\Lambda_i^j)\prod_{\beta=1}^3 |E_i|^{1/p_\beta},$$
where $\tilde{C}(\Lambda_i^j)$ are defined in (\ref{C}) and (\ref{C2}).
So $\Lambda_i^j$ is of restricted weak type and we have the following estimate about $C(\Lambda_i^j)$~:
$$ C(\Lambda_i^j)\lesssim \tilde{C}(\Lambda_i^j). $$
In addition the parameter $N$ in (\ref{C}) and (\ref{C2}) is as large as we want, and we have also obtained the desired estimates on $C(\Lambda_i^j)$.
\findem

\mb The proof of Theorem \ref{thm:res2} is now completed. By using the concept of ``restricted weak type'', we can have
a ``stronger'' result than Theorem \ref{thm:central}~:

\begin{thm} \label{thfort}  Let $T,p,q,r$ be an operator and exponents of Theorem \ref{thm:central}.
Then for all $\delta\geq 1$, there exist a constant $C=C(p,q,r,\delta)$ (independent on the interval $I$) such that for
all sets $E_3$ of finite measure, there exists a substantial subset $E'_3 \subset E_3$ satisfying that for all
functions $f\in \s(\R),\ g\in \s(\R)$ and $h\in F(E_3')$ we have
\begin{align*}
\lefteqn{\left|\langle T(f,g),h{\bf 1}_{I} \rangle \right| \leq} & & \\
 & &  C \left(\sum_{k\geq 0} 2^{-k(1/p+\delta)} \|f{\bf 1}_{2^kI}\|_p\right) \left(\sum_{k\geq 0} 2^{-k(1/q+\delta)} \|g{\bf 1}_{2^kI}\|_q\right)  |E_3|^{1/r'}.
\end{align*}
\end{thm}

\mb When $r>1$, this result is stronger than Theorem \ref{thm:central} but less practicable. We now prove it because it
will be useful in the sequel.

\dem The proof is exactly the same as the previous one, so we shall only explain the modifications. We always study the
trilinear form
$$ \Lambda(f,g,h) := \langle T(f,g) , h {\bf 1}_I \rangle. $$
In the section \ref{discretis} we have seen that the study of $\Lambda$ can be reduced to the study of the following model sum
$$ \Lambda(f,g,h) =  \sum_{s\in \mathbf{S}} |I_s|^{-1/2}\epsilon_{s} \langle
\phi_{s_1},f \rangle \langle \phi_{s_2},g\rangle \langle \phi_{s_3},h {\bf 1}_{I} \rangle, $$
where $\mathbf{S}$ is a general collection of tri-tiles.
Then we have decomposed this sum with (\ref{decompop}) by~:
$$\Lambda(f,g,h) =  \sum_{k_1,k_2\geq 0} \Lambda_0^{k_1,k_2}(f,g,h) + \sum_{\genfrac{}{}{0pt}{}{k_1,k_2\geq 0}{l\leq 0}} \Lambda_1
^{k_1,k_2,l}(f,g,h).$$ By Theorem \ref{thm:res}, we have shown that the trilinear forms $\Lambda_i^j$ are of restricted
weak type $(p,q,r')$ and we have obtained estimates on their bounds. The construction of the substantial subset
$E_\alpha'=E_3'$ does not depend on the trilinear form $\Lambda_i^j$, so we can deduce that our trilinear form
$\Lambda$ is always of restricted weak type. Also for measurable sets $E_1,E_2,E_3$ of finite measure, there exists a
substantial subset $E_3' \subset E_3$ such that for all functions $f\in F(E_1),g\in F(E_2)$ and $h\in F(E_3')$ we
have~:
$$ \left| \Lambda(f,g,h) \right| \lesssim  |E_3|^{1/r'} \left[\sum_{k_1,k_2\geq 0} 2^{-\delta'(k_1+k_2)} \left|E_1 \cap C_{k_1}\right|^{1/p} \left|E_2 \cap C_{k_2}\right|^{1/q} \right].$$
Here $\delta'$ is an exponent as large as we want. Over each corona, by using the real interpolation on the exponents
$p$ and $q$ (so $r$ is fixed), we obtain also the desired result. \findem

\mb Having obtained our main result for the $x$-independent symbols, we will extend our result for maximal operators
and for $x$-dependent symbols in the next section.

\section{More general bilinear operators.}
\label{sectiont}

\mb Let us name our ``off-diagonal'' estimates for convenience.

\begin{df} Let $T$ be an operator (maybe non bilinear) acting from $\s(\R) \times \s(\R)$ into $\s'(\R)$. For $0<p,q,r\leq \infty$ exponents such that
$$ \frac{1}{r}= \frac{1}{p}+\frac{1}{q},$$
we say that $T$ satisfies {\it ``off-diagonal'' estimates at the scale $L$ and at the order $\delta$}, in short $T\in
{\mathcal O}_{L,\delta}(L^p \times L^q,L^r)$, if there exists a constant $C=C(p,q,r,L,\delta)$ such that for all
functions $f,g\in\s(\R)$ and all interval $I$ of length $|I|=L$, we have
\be{offdiag} \|T(f,g)\|_{r,I} \leq C \left[
\sum_{k\geq 0} 2^{-k(\delta+1/p)} \|f\|_{p,2^{k+1}I} \right]\left[ \sum_{k\geq 0} 2^{-k(\delta+1/q)} \|g\|_{q,2^{k+1}I}
\right]. \ee
\end{df}

\begin{rem} Equivalently, an operator $T$ satisfies ``off-diagonal'' at the scale $L$ and at the order $\delta$ if there exists a constant $C=C(p,q,r,L,\delta)$ with for all functions $f,g\in\s(\R)$ and all interval $I$ of length $|I|=L$, we have
$$ \|T(f,g)\|_{r,I} \leq C \left[ \sum_{k\geq 0} 2^{-k(\delta+1/p)} \|f\|_{p,C_k(I)} \right]\left[ \sum_{k\geq 0} 2^{-k(\delta+1/q)} \|g\|_{q,C_k(I)} \right]. $$
This is a better way to describe the ``off-diagonal'' decay of an operator $T$ and these properties can been describe
as in Corollary \ref{corrolaire}.
\end{rem}

\mb First we generalize the previous result for maximal operators.

\subsection{``Off-diagonal'' estimates for maximal bilinear operators.}

We have the three following theorems~:

\begin{thm} \label{max} Let $\Delta$ be a nondegenerate line in the frequency plane. Let $1<p,q \leq\infty$ be exponents such that
$$0<\frac{1}{r}=\frac{1}{q}+\frac{1}{p}<\frac{3}{2}.$$
For all $\delta\geq 1$, $L>0$, for all symbol $\sigma$ supported in $\{(\alpha,\beta),\ d((\alpha,\beta),\Delta)\geq L^{-1}\}$ satisfying
$$ \forall b,c\geq 0 \qquad \left|\partial^{b}_\alpha \partial^c_\beta \sigma(\alpha,\beta)\right| \lesssim |d((\alpha,\beta),\Delta)|^{-b-c} $$
and for all smooth function $\phi$, which is equal to $1$ around $0$, the maximal bilinear operator \be{opbi2}
T_{max}(f,g)(x):=\sup_{r>0} \left|\int e^{ix(\alpha+\beta)} \widehat{f}(\alpha)\widehat{g}(\beta)\sigma(\alpha,\beta)
\Big[1-\phi(r(\alpha-\beta))\Big] d\alpha d\beta\right|\ee satisfies ``off-diagonal'' estimates at the scale $L$ and at
the order $\delta$ : $T_{max}\in {\mathcal O}_{L,\delta}(L^p \times L^q,L^r)$. In addition the implicit constant can be uniformly bounded in $L>0$.
\end{thm}

\begin{thm} \label{max1}
For the same exponents, we have the same continuities for the maximal bilinear operator (at the scale $L$)~:
$$M^{L}(f,g)(x):=\sup_{0<r\leq L}\, \frac{1}{r}\int_{|t|\leq r} \left|f(x-t) g(x+t) \right| dt.$$
\end{thm}

\begin{thm} \label{max2}
Let $K$ be a kernel on $\R$ satisfying the H\"ormander's conditions, then the maximal bilinear operator \be{opbi3}
T^{L}_{max}(f,g)(x):=\sup_{0<\epsilon<r<L} \left|\int_{\epsilon\leq |y|\leq r} f(x-y)g(x+y) K(y) dy \right|  \ee
verifies the same local estimates : $T^L_{max}\in {\mathcal O}_{L,\delta}(L^p \times L^q,L^r)$ (for the exponents $p,q,r$ of Theorem \ref{max}).
\end{thm}

\dem The proof of these three Theorem is a shake between the proof of our Theorem \ref{thm:central} and an additional maximality argument.
The maximal truncation in the physical space (Theorems \ref{max1} and \ref{max2}) is a little more complex than the maximal truncation in the frequency space (Theorem \ref{max}). So we deal with the two last theorems and just explain the modifications to prove them. 
The maximal version of the different arguments
have first shown in \cite{lacey} by M. Lacey and then improved in \cite{dtt} by C. Demeter, T. Tao and C. Thiele. In these articles, the authors study the behaviour of the maximal averages (like Theorem \ref{max1}). The remark 1.6 of \cite{dtt} specifies the similarity between the operators of Theorem \ref{max1} and \ref{max2}.
So in fact our three previous theorems are an illustration of the same ideas, we will not detail. \\
The subsection \ref{discretis} is based on the decomposition of the bilinear operator by discrete models. For our maximal operators, the same reduction is shown in \cite{dtt} (Theorem 4.4) and the important condition (\ref{remarque}) for the tiles are
always satisfied. Then the maximal version of Proposition \ref{prop} is written in \cite{dtt} too (there is a new
factor in the different inequalities but it is not important). We have exactly the same version of Lemma \ref{arbre}
for maximal bilinear operators (see Proposition 6.2 of \cite{dtt}). Using these few technical modifications, we can
compute our proof of Theorem \ref{thm:central} and obtain its maximal versions : Theorems \ref{max},\ref{max1} and \ref{max2}. \findem

\subsection{Proof of Theorem \ref{thm:central} for $x$-dependent symbols.}

\mb In this subsection, we prove the ``off-diagonal'' estimates of Theorem \ref{thm:central} in the case where the
symbol $\sigma$ depends on the spatial variable $x$ and also we complete the proof of our main result.

\begin{thm} Let $\Delta$ be an nondegenerate line of the frequencial space. Let $\sigma\in C^\infty(\R^3)$ be a symbol satisfying~:
$\forall a,b,c\geq 0$
$$ \left| \partial_x^a \partial_\alpha ^b \partial_\beta^c \sigma(x,\alpha,\beta) \right| \lesssim \left(1+d((\alpha,\beta),\Delta) \right)^{-b-c}.$$
Then the bilinear operator $T_\sigma$ (defined on $\s(\R)\times \s(\R)$ by (\ref{operateur})) verifies $T\in {\mathcal
O}_{1,\delta}(L^p \times L^q,L^r)$ for any $\delta\geq 0$ and any exponents $p,q,r$ such that
$$0<\frac{1}{r}=\frac{1}{p}+\frac{1}{q} < \frac{3}{2} \textrm{   and   } 1<p,q\leq \infty.$$
\end{thm}

\mb Our assumptions for the symbol corresponds to the class $BS^{1,0}_{\theta}$ of \cite{bipseudo} where the angle
$\theta\in]-\pi/2,\pi/2[\setminus\{0,-\pi/4\}$ is given by the line $\Delta$~: $$ \Delta:=\left\{ (\alpha,\beta),\
\beta=\tan(\theta) \alpha \right\}.$$ For convenience, we will deal in the proof only with the case $\theta=\pi/4$. The
important fact is that the singular quantity $\beta-\tan(\theta)\alpha$ does not correspond to the quantity
$\alpha+\beta$, which appears in the exponential term of (\ref{operateur}). The limit and particular case
$\theta=-\pi/4$ is studied in \cite{bipseudo}.

\dem \\ The proof is quite technical. We will also assume that $r\geq 1$ (which allows us to simplify a few arguments). Then we will explain in Remark \ref{modifs} how modify the proof to obtain the same result when $r<1$. \\
So we fix an interval $I$ of length $|I|=1$. We use a decomposition of the symbol $\sigma$. Let $\Phi$ be a smooth function on $\R$ such that
$$ |x|\leq 1 \Longrightarrow \phi(x)=1 \textrm{  and  supp} (\phi) \subset [-2,2].$$
We have also
\begin{eqnarray*}
 \sigma(x,\alpha,\beta) & = & \sigma(x,\alpha,\beta) [1-\Phi(\alpha-\beta)]+ \sigma(x,\alpha,\beta) \Phi(\alpha,\beta) \\
 & := &  \sigma^{\infty}(x,\alpha,\beta) + \sigma^0(x,\alpha,\beta).
\end{eqnarray*}
$1-)$ The case of the symbol $\sigma^\infty$. \\
We have an operator associated to this symbol~:
\begin{align*}
T^\infty(f,g)(x):=  \int_{\R^2} e^{ix(\alpha+\beta)}
\widehat{f}(\alpha)\widehat{g}(\beta)\sigma(x,\alpha,\beta)[1-\Phi(\alpha-\beta)] d\alpha d\beta,
 \end{align*}
which can been written
$$T^\infty(f,g)(x)= U_{x}(f,g)(x),$$
with $U$ defined by
\begin{align*}
U_{y}(f,g)(x):= \int_{\R^2} e^{ix(\alpha+\beta)} \widehat{f}(\alpha)\widehat{g}(\beta)\sigma(y,\alpha,\beta)[1-\Phi(\alpha-\beta)] d\alpha d\beta.
\end{align*}
By using the Sobolev's imbedding $W^{1,r}(I) \hookrightarrow L^\infty(I)$ (because $r\geq 1$), we get
$$ \forall x\in I,\quad \left| T^\infty(f,g)(x) \right| \leq \| U_{y}(f,g)(x)\|_{\infty,y\in I} \lesssim \sum_{k=0}^{1} \| \partial_y^k U_{y}(f,g)(x) {\bf 1}_{I}(y)\|_{r,dy}.$$
Then by integrating for $x\in I$ and using Fubini's Theorem, we obtain
$$ \left\| T^\infty(f,g)\right\|_{r,I} \lesssim \sum_{k=0}^{1}  \left\| \left\|\partial_y^k U_{y}(f,g) \right\|_{r,I} \right\|_{r,I,dy}.$$
We can fix $k\in \{0,1\}$ and $y\in I$. Then we have
$$\left\|\partial_y^k U_{y}(f,g) \right\|_{r,dx} \lesssim \left\| V(f,g) \right\|_{r,I}, $$
where $V$ is the bilinear operator defined by
$$V(f,g)(x):= \int_{\R^2} e^{ix(\alpha+\beta)} \widehat{f}(\alpha)\widehat{g}(\beta) \partial_{y}^{k} \sigma(y,\alpha,\beta)[1-\Phi(\alpha-\beta)] d\alpha d\beta.$$
So $V=T_\tau$ is the bilinear operator associated to the $x$-independent symbol
$$\tau(\alpha,\beta):= \partial_{y}^{k} \sigma(y,\alpha,\beta)[1-\Phi(\alpha-\beta)].$$
From the assumptions about $\sigma$, the symbol $\tau$ satisfies : for all $b,c\geq 0$
$$\left|\partial_\alpha^b \partial_\beta^c \tau(\alpha,\beta) \right|\lesssim \left|\alpha-\beta\right|^{-n-p} .$$
In addition $\tau$ is supported in the domain~: $\{|\alpha-\beta|\geq 1\}$. We can also apply Theorem \ref{thm:central}
proved in section \ref{section1} for $x$-independent symbol. For all $\delta\geq 1$, we have an ``off-diagonal''
estimate at the scale $1$~:
$$\left\| V(f,g)\right\|_{r,I} \lesssim \left(\sum_{k_1\geq 0} 2^{-k_1(1/p+\delta)} \|f\|_{p,2^{k_1}I}\right)
\left(\sum_{k_2\geq 0} 2^{-k_2(1/q+\delta)} \|g\|_{q,2^{k_2}I}\right).$$ All theses estimates are uniform with respect to
$k\in\{0,1\}$ and $y\in I$, so we get \be{estitruc}  \left\| T^\infty(f,g)\right\|_{r,I} \lesssim  \left(\sum_{k_1\geq
0} 2^{-k_1(1/p+\delta)} \|f\|_{p,2^{k_1}I}\right) \left(\sum_{k_2\geq 0} 2^{-k_2(1/q+\delta)} \|g\|_{q,2^{k_2}I}\right)
.\ee
So we have shown the desired estimates for this first term.\\
$2-)$ The case of the symbol $\sigma^0$. \\
The associated operator is given by
$$T^0(f,g)(x):=\int_{\R^2} e^{ix(\alpha+\beta)} \widehat{f}(\alpha)\widehat{g}(\beta)\sigma(x,\alpha,\beta) \Phi(\alpha,\beta) d\alpha d\beta.$$
We use the same arguments as for the first point. So we have to study the operator $V$ defined by
$$V(f,g)(x):= \int_{\R^2} e^{ix(\alpha+\beta)} \widehat{f}(\alpha)\widehat{g}(\beta) \partial_{y}^{k} \sigma(y,\alpha,\beta)\Phi(\alpha-\beta)d\alpha d\beta.$$
The parameters $k\in\{0,1\}$ and $y\in I$ are fixed. The symbol associated to this operator is supported on $\{(\alpha,\beta),\ |\alpha-\beta|\leq 2\}$. That is why, we use modulations to move this support.
\begin{align*}
V(f,g)(x) & =  \int_{\R^2} e^{ix(\alpha+\beta)} \widehat{f}(\alpha +3)\widehat{g}(\beta-3) \partial_{y}^{k} \sigma(y,\alpha+3,\beta-3)\Phi(\alpha-\beta +6)  d\alpha d\beta \\
 & = \int_{\R^2} e^{ix(\alpha+\beta)} \widehat{e^{3i.} f}(\alpha)\widehat{e^{-3i.} g}(\beta) \partial_{y}^{k} \sigma(y,\alpha+3,\beta-3)\Phi(\alpha-\beta +6) d\alpha d\beta.
\end{align*}
Also $V$ is now the bilinear operator, applied to the functions $e^{3i.} f$ and $e^{-3i .} g$, whose 
($x$-independent) symbol
$$\tau(\alpha,\beta):= \partial_{y}^{k} \sigma(y,\alpha+3,\beta-3)\Phi(\alpha-\beta+6).$$
is supported on
$$\{(\alpha,\beta),\ |\alpha-\beta+6|\leq 2 \} \subset \{(\alpha,\beta),\ 1\leq |\alpha-\beta| \leq 8\}$$
and satisfies
\begin{align*}
 \forall b,c\geq 0, \qquad \left|\partial^{b}_\alpha \partial^c_\beta \tau(\alpha,\beta)\right| &  \lesssim \max_{0\leq j\leq b} \max_{0\leq i\leq c}  \left(1+|\alpha-\beta+6|\right)^{-i-j} {\bf 1}_{1\leq |\alpha-\beta| \leq 8}  \\
 & \lesssim {\bf 1}_{1\leq |\alpha-\beta| \leq 8} \\
 & \lesssim {\bf 1}_{1\leq |\alpha-\beta| \leq 8} \left|\alpha-\beta\right|^{-b-c}.
\end{align*}
Also we can again use Theorem \ref{thm:central} (proved in Section \ref{section1} for $x$-independent symbol) and we
obtain~:
$$\left\| V (f,g)\right\|_{r,I} \lesssim \left(\sum_{k_1\geq 0} 2^{-k_1(1/p+\delta)} \| f\|_{p,I}\right)
\left(\sum_{k_2\geq 0} 2^{-k_2(1/q+\delta)} \|g\|_{q,I}\right).$$
We have also finished the proof.
\findem

\mb
\begin{rem} \label{modifs}
In this remark, we want to explain how modify the previous proof when $r<1$. When we want to study bilinear operators
with $r<1$, we have to use the associated trilinear form and the concept of "restricted weak type" (defined in
Definition \ref{restricteddef}). These two arguments allow us to get round the lack of the triangular inequality in the
space $L^{r}$. Let
$$\Lambda(f,g,h):= \langle T(f,g),h\rangle.$$
We have
$$ \Lambda(f,g,h) = \int_\R \int_{\R^2} e^{ix(\alpha+\beta)} \sigma(x,\alpha,\beta) \widehat{f}(\alpha)
\widehat{g}(\beta) h(x) d\alpha d\beta dx.$$ We use the same decomposition of $\sigma$, getting the trilinear forms
$\Lambda^\infty$ and $\Lambda^0$. Let us study first $\Lambda^\infty$ and fix an interval $I$ of
length $|I|=1$. We take a function $h\in \s(\R)$, which is supported on $I$. We use again the Sobolev imbedding
$W^{1,1}(I) \hookrightarrow L^\infty(I)$. By writing
$$ \left|\Lambda^\infty(f,g,h) \right| \leq \int_\R \left\| U_{y}(f,g)(x) {\bf 1}_{I}(y)\right\|_{\infty,I,dy} \left|h(x) \right|{\bf 1}_{I}(x) dx.$$
We can also obtain
$$ \left|\Lambda^\infty(f,g,h) \right| \lesssim \int_I \int_I \left| U_{y}(f,g)(x)\right|\left|h(x) \right|dx dy + \int_I \int_I \left| \partial_y U_{y}(f,g)(x)\right|\left|h(x) \right|dx dy .$$
Then when $y\in I$ and $k\in\{0,1\}$ are fixed, we find again the quantities
$$ \int_I \left| \partial_y^kU_{y}(f,g)(x)\right|\left|h(x) \right|dx .$$
Now the bilinear operator $\partial_y^k U_{y}$ is associated to an $x$-independent symbol, which verifies the good
assumptions. We can also use Theorem \ref{thfort} in order to obtain the wished estimates (\ref{estitruc}) in a
``restricted weak type sense'' for the exponent $r$. We produce the same modifications to study $\Lambda^0$. By
noticing that the way to construct the substantial subset (in the definition of restricted weak type) does not depend
on the trilinear form, we can deduce that the trilinear form $\Lambda$ satisfies (\ref{estitruc}) in a ``restricted
weak type sense'' too. Then we use interpolation on the exponent $r$, to obtain exactly (\ref{estitruc}), which allows
us to conclude.
\end{rem}

\section{Continuities for bilinear operators, satisfying ``off-diagonal'' estimates.}
\label{weight}

\mb Recall that in the linear case, by using the maximal sharp function, we can prove weighted continuities for
linear operator with the Muckenhoupt's weights. In the bilinear case, we do not have a good substitute to the maximal
sharp function. That is why we shall use the previous ``off-diagonal'' estimates to obtain weighted global continuities
on Lebesgue spaces and in particular prove Theorem  \ref{thmcentral}.

\mb We first want to give an application of these off-diagonal estimates. Recall that in the previous sections, we
have proved that our bilinear operators (and maximal bilinear operators) satisfy these ``off-diagonal'' estimates at
any order. The time-frequency analysis does not work for functions in the $L^\infty$ space. So we don't know if
our operators $T$ are bounded from $L^\infty \times L^\infty$ in $BMO$. However these local estimates give a weak
result about the behavior of $T(f,g)$ when the two functions $f$ and $g$ belong to $L^\infty$.

\begin{prop} Let $f,g$ be two functions of $L^1(\R) \cap L^\infty(\R)$ and fix $1<r<\infty$. If there exist $L>0,\delta\geq 1$ and $p,q>1$ such that an operator $T \in {\mathcal O}_{\delta,L}(L^p\times L^q,L^r)$, then we have
$$ \limsup_{|I| \rightarrow \infty}  \left(\frac{1}{|I|}\int_{I} | T (f,g)|^{r} \right)^{1/r} \lesssim \|f\|_\infty \|g\|_\infty.$$
Here we take the limit when $I$ is an interval with $|I|\to \infty$ and the implicit constant does not depend on the
two functions $f$ and $g$ and on the parameter $L$.
\end{prop}

\dem We set $I_i:=[iL,(i+1)L[$ for all $i\in \Z$. Then for $I$ with $|I|>>L$, we get
$$ \int_{I} | T(f,g)|^{r} \leq \sum_{\genfrac{}{}{0pt}{}{i\in\Z}{I_i\cap I \neq \emptyset}} \int_{I_i} | T(f,g)|^{r}.$$
However the number of index $i$ which appears in the sum is bounded by $|I|/L$, so by using the local estimate we get
\begin{align*}
\int_{I} | T(f,g)|^{r} & \lesssim \sum_{\genfrac{}{}{0pt}{}{i\in\Z}{I_i\cap I\neq \emptyset}} \frac{L}{|I_i|}\int_{I_i} | T(f,g)|^{r} \\
 & \lesssim \sum_{\genfrac{}{}{0pt}{}{i\in\Z}{I_i\cap I \neq \emptyset}} L \|f\|_\infty^{r} \|g\|_\infty^{r} \\
 & \lesssim |I| \|f\|_\infty^{r} \|g\|_\infty^{r}.
\end{align*}
The second inequality is due to the fact that~:
$$ |I_i|^{1/r}\|T(f,g)\|_{r,I_i} \lesssim \left[\inf_{x\in I_i}
M_{HL}(f)(x) \right]\left[\inf_{x\in I_i} M_{HL}(g)(x) \right] \lesssim \|f\|_\infty \|g\|_\infty.$$
So we obtain
$$ \left(\frac{1}{|I|}\int_{I} | T_{max} (f,g)|^{r} \right)^{1/r} \lesssim \|f\|_\infty \|g\|_\infty, $$
uniformly with $L$ for $|I|$ large enough.
\findem

\mb Let us now define our weights.
\begin{df} Let $\theta>0$ and $l>0$ be fixed. We set that a nonnegative function $\omega$ belongs to the {\it class $\po_\theta(l)$} if there exists a constant $C$ such that for all interval $I$ of length $|I|=l$ and for all integer $k\geq 0$, we have
\be{poids} 2^{-k\theta} \sup_{x\in I} \omega(x) \leq C \inf_{2^{k}I} \omega(x).\ee
\end{df}

\mb So we claim that a function $\omega\in \po_{\theta}(l)$ seems to be likely a polynomial function whose the degree
is less than $\theta$ and it is almost constant at the scale $l$. We show in the following example that these classes
are not empty.

\begin{ex} For all $\theta>0$ and $\alpha\in[0,\theta[$, the functions
$$x\rightarrow 1, \qquad x\rightarrow (1+|x|)^{\alpha} \qquad \textrm{ and } \qquad x\rightarrow (1+|x|)^{-\alpha}$$
belong to the class $\po_\theta(1)$. The proof is easy and left to the reader.
\end{ex}

\begin{rem} \label{croissancepol} In fact, it is easy to prove that a weight $\omega$ belongs to the class $\po_\theta (l)$ if and only if there exists a constant $C$ such that
$$ \forall x,y\in\R, \qquad \omega(x) \leq C \left( 1+ \frac{|x-y|}{l} \right)^{\theta} \omega(y).$$
We cannot compare these weights with the Muckenhoupt's weights, because for $\omega \in \po_\theta(l)$ we have informations only at the scale $l$.
\end{rem}

\mb We have the following result~:

\begin{thm} \label{thpoids} Let $T$ be a bilinear operator and $0<p,q,r<\infty$ be exponents satisfying
$$\frac{1}{r}=\frac{1}{p}+\frac{1}{q} \qquad \textrm{and} \qquad 1\leq p,q.$$
For $\delta>0$ and $l>0$, if  $T$ satisfies ``off-diagonal'' estimates at the order $\delta$ and at the scale $l$, then
for all $\omega\in \po_\theta(l)$ with $0\leq \theta<\delta \max\{r,1\}$, the operator $T$ is continuous from
$L^{p}(\omega) \times L^q(\omega)$ in $L^{r}(\omega)$.
\end{thm}

\dem To check this, recall that for all interval $I$ of length $|I|=l$, \be{hyppoids} \left(\int_I \left| T(f,g)
\right|^r\right)^{1/r} \lesssim \left(\sum_{k\geq 0} 2^{-k(1/p+\delta)} \|f\|_{p,2^{k}I}\right) \left(\sum_{k\geq 0}
2^{-k(1/q+\delta)} \|g\|_{q,2^{k} I}\right).\ee So we decompose the whole space $\R$ with the disjoint intervals
$I_i$ defined by $I_i=[il,(i+1)l[$ for $i\in \Z$. So we have
$$\left\| T(f,g)\right\|_{r,wdx} = \left\| \left\| T(f,g)\right\|_{r,wdx,I_i} \right\|_{r,i\in\Z}.$$
Let $i\in \Z$ be fixed, we use (\ref{poids}) and (\ref{hyppoids}) to obtain
\begin{align*}
\left\| T(f,g)\right\|_{r,wdx,I_i} & \leq \left\| w\right\|_{\infty,I_i}^{1/r} \left\| T(f,g) \right\|_{r,I_i} \\
 & \lesssim  \left\| w\right\|_{\infty,I_i}^{1/r} \left(\sum_{k\geq 0} 2^{-k(1/p+\delta)} \|f\|_{p,2^{k}I_i}\right) \left(\sum_{k\geq 0} 2^{-k(1/q+\delta)} \|g\|_{q,2^{k} I_i}\right)
 \end{align*}
We estimate the first sum with
\begin{align*}
\left\| w\right\|_{\infty,I_i}^{1/p} \left(\sum_{k\geq 0} 2^{-k(1/p+\delta)} \|f\|_{p,2^{k}I_i}\right) & \lesssim   \left(\sum_{k\geq 0} 2^{-k(1/p+\delta)} \left\| w\right\|_{\infty,I_i}^{1/p}\|f\|_{p,2^{k}I_i}\right) \\
 & \lesssim   \left(\sum_{k\geq 0} 2^{-k(1/p+\delta)} 2^{k\theta/p} \inf_{2^{k}I_i} \omega ^{1/p}\|f\|_{p,2^{k}I_i}\right)  \\
 & \lesssim   \left(\sum_{k\geq 0} 2^{-k(1/p+\delta-\theta/p)} \|f\|_{p,wdx,2^{k}I_i}\right).
\end{align*}
The second term is studied by the same way. By summing over $i\in\Z$, we get
\begin{align*}
\lefteqn{\left\| T(f,g)\right\|_{r,wdx} \lesssim} & & \\
 & &  \left\| \left(\sum_{k\geq 0} 2^{-k(1/p+\delta-\theta/p)} \|f\|_{p,wdx,2^{k}I_i}\right) \left(\sum_{k\geq 0} 2^{-k(1/q+\delta-\theta/q)} \|g\|_{q,wdx,2^{k} I_i}\right)  \right\|_{r,i\in\Z}.
\end{align*}
With the help of Hölder's and Minkowski's inequalities, we get
\begin{align*}
\left\| T(f,g)\right\|_{r,wdx} \lesssim  &  \left(\sum_{k\geq 0} 2^{-k(1/p+\delta-\theta/p)} \left\| \|f\|_{p,wdx,2^{k}I_i}\right\|_{p,i\in\Z}\right) \\
 & \left(\sum_{k\geq 0} 2^{-k(1/q+\delta-\theta/q)} \left\| \|g\|_{q,wdx,2^{k} I_i} \right\|_{q,i\in\Z}\right).
\end{align*}
However the collection of sets $(2^{k}I_i)_i$ is a $2^{k}$-covering, so
\begin{align*}
\left\| T(f,g)\right\|_{r,wdx} & \lesssim \left(\sum_{k\geq 0} 2^{-k(\delta-\theta/p)}  \|f\|_{p,wdx}\right) \left(\sum_{k\geq 0} 2^{-k(\delta-\theta/q)}  \|g\|_{q,wdx}\right).
\end{align*}
Then we conclude with the fact that $p,q>1$ hence
$$\max\left\{\frac{\theta}{p},\frac{\theta}{q}\right\}  \leq \left\{\begin{array}{ll}
   \frac{\theta}{r} <\delta & \textrm{ if $r\geq 1$} \\
   \theta <\delta & \textrm{ if $r\leq 1$}
  \end{array}\right. .$$
\findem

\begin{rem} From the fact that the weight $\omega(x)=1$ belongs to the class $\po_\theta(L)$, we have also proved that the operators of Theorem \ref{thm:central} and the maximal operators of Theorems \ref{max},\ref{max1}, \ref{max2} are bounded in classical Lebesgue spaces.
\end{rem}

\mb We complete this result by the following proposition in Sobolev spaces.

\begin{df} Let $\omega$ be a weight on $\R$. For all $m\geq 0$ and $p\in]1,\infty[$, we set $W^{m,p}(\omega)$ for the Sobolev space on $\R$ with the weight $\omega$, defined as the set of distributions $f\in\s'(\R)$ such that $J_m(f) \in L^{p}(\omega)$, where $J_m:=\left(Id-\Delta \right)^{m/2}$.
\end{df}

\begin{prop} Let $\Delta$ be a nondegenerate line, $\omega$ be a weight in $\cup_{\theta\geq 0} \po_\theta(1)$ and $\sigma\in C^\infty(\R^3)$ be a symbol satisfying
$$ \forall a,b,c\geq 0, \qquad  \left| \partial_x^a \partial_\alpha ^b \partial_\beta^c \sigma(x,\alpha,\beta) \right| \lesssim \left(1+d((\alpha,\beta),\Delta) \right)^{-b-c}.$$
Let $p,q,r$ be exponents satisfying
$$ 0<\frac{1}{r}=\frac{1}{p}+\frac{1}{q}<\frac{3}{2} \quad \textrm{and} \quad 1<p,q\leq \infty.$$
Then the bilinear operator $T_\sigma$ (defined on $\s(\R)\times \s(\R)$ by (\ref{operateur})) satisfies~:
 for all integer $n\geq 0$
\be{derivation} \forall f,g\in\s(\R),\qquad \left\| D^{(n)} T_\sigma(f,g) \right\|_{L^{r}(\omega)} \lesssim \sum_{\genfrac{}{}{0pt}{}{0\leq i,j \leq n}{i+j\leq n}} \|D^{(i)} f\|_{L^p(\omega)} \|D^{(j)} g\|_{L^q(\omega)}. \ee
 Here we write $D^{(i)}$ for the derivative operator of order $i$. Also $T_\sigma$ is continuous from $W^{m,p}(\omega) \times W^{m,q}(\omega)$ in $W^{m,r}(\omega)$ for all real $m\geq 0$.
\end{prop}

\dem Let us begin to prove (\ref{derivation}). The two functions $f$ and $g$ are smooth so we can differentiate the integral defining $T_\sigma(f,g)$. It is also easy to check that 
$$D^{(1)}T_\sigma(f,g)= T_\sigma(D^{(1)}f,g) + T_\sigma(f,D^{(1)}g) + T_{\partial_x\sigma}(f,g).$$
Then for higher orders, we get
$$D^{(n)}T_\sigma(f,g)= \sum_{\genfrac{}{}{0pt}{}{0\leq i,j,k \leq n}{i+j+k=n}} T_{\partial_x^k\sigma}(D^{(i)}f,D^{(j)}g).$$
By using the previous Theorems \ref{thm:central} and \ref{thpoids}, we obtain (\ref{derivation}). We can also deduce a weaker estimate~:
$$ \forall f,g\in\s(\R),\qquad \left\| D^{(n)} T_\sigma(f,g) \right\|_{r,\omega} \lesssim  \|f\|_{W^{n,p}(\omega)} \|g\|_{W^{n,q}(\omega)}.$$
By density (see Lemma \ref{lemmesobolev}), the operator $T_\sigma$ can be continuously extended from $W^{n,p}(\omega)\times W^{n,q}(\omega)$ into
$W^{n,r}(\omega)$. Then we will use interpolation to extend this result when $n$ is not an integer. The exponents
$p,q,r$ are fixed and we study the bilinear operator $T_\sigma$. We have shown that $T_\sigma$ is continuous from
$W^{n,p}(\omega) \times W^{n,q}(\omega)$ into $W^{n,r}(\omega)$, for all integer $n$. By using the bilinear
interpolation (with Lemme \ref{lemmesobolev}) \footnote{The theory of the multilinear interpolation is studied in the chapter $4$ of \cite{LP} for the
real interpolation and in Theorem 4.4.1 of \cite{BL} for the complex interpolation.} on $n$, we finish the proof.
\findem

\begin{lem} \label{lemmesobolev} For all weight $\omega\in \cup_{\theta\geq 0} \po_\theta(1)$, all exponent $1<p<\infty$ and all real $s\geq 0$, the space $\s(\R)$ is a dense subspace in $W^{s,p}(\omega)$. In addition, the collection of Sobolev spaces $(W^{s,p}(\omega))_{s\geq 0}$ form an interpolation scale.
\end{lem}

\dem Let $\omega$ be a fixed weight in $\cup_{\theta\geq 0} \po_\theta(1)$. We have seen in Remark \ref{croissancepol}
that $\omega$ has a polynomial growth. Due to the fact that $J_s(\s(\R)) = \s(\R)$, we have the inclusion $\s(\R) \subset W^{s,p}(\omega)$. Recall that $J_s:=\left(Id-\Delta \right)^{s/2}$. In addition, we have that $L^p(\omega) \subset \s'(\R)$, so we can compute the operator $J_{-s}$ on the space $L^p(\omega)$. We finally obtain that $J_s$ is an automorphism from $W^{s,p}(\omega)$ to $L^p(\omega)$ and an isomorphism on $\s(\R)$. As $\s(\R)$ is dense in $L^p(\omega)$, we get the density of $\s(\R)$ into the Sobolev space $W^{s,p}(\omega)$. \\
For the interpolation claim, we omit the details. The classical proof for complex interpolation with $\omega=1$ can easily been extended to the general case. 
\findem

\mb

\begin{rem} From the fact that the weight $\omega(x)=1$ belongs to the class $\po_\theta(1)$, we have also proved that the operators of Theorem \ref{thmcentral} satisfy an Hölder's inequality in Sobolev spaces. 
\end{rem}

\mb

\begin{rem}
Also with the notation of \cite{bipseudo}, we have proved continuities for all the operators associated to symbols
$\sigma\in BS_{1,0;\theta}^{0}$. In addition, we have described the action of these operators on Sobolev spaces.
 This is an interesting improvement of this article and it incite us to obtain new
results in order to continue the construction of a bilinear pseudodifferential calculus. We will do it in a next paper
\cite{pseudo} by introducing new larger symbolic classes of bilinear symbols of order $(m_1,m_2)$.\\
About continuities in Lebesgue spaces, a question is still open : what about the class $BS^0_{\rho,\delta;\theta}$
(defined in \cite{bipseudo}) ?
\end{rem}

\end{document}